# A GENERAL TRIMMING APPROACH TO ROBUST CLUSTER ANALYSIS[1]

By Luis A. García-Escudero, Alfonso Gordaliza,
Carlos Matrán and Agustin Mayo-Iscar

*Universidad de Valladolid*

We introduce a new method for performing clustering with the aim of fitting clusters with different scatters and weights. It is designed by allowing to handle a proportion $\alpha$ of contaminating data to guarantee the robustness of the method. As a characteristic feature, restrictions on the ratio between the maximum and the minimum eigenvalues of the groups scatter matrices are introduced. This makes the problem to be well defined and guarantees the consistency of the sample solutions to the population ones.

The method covers a wide range of clustering approaches depending on the strength of the chosen restrictions. Our proposal includes an algorithm for approximately solving the sample problem.

**1. Introduction.** Many statistical practitioners view cluster analysis as a collection of mostly heuristic techniques for partitioning multivariate data. This arises from the fact that most cluster techniques are not *explicitly* based on a probabilistic model, and could lead to the feeling that no assumption is necessary and that the obtained results are "objective" (see the comments on page 123 in Flury [7]). However, objectiveness is far from reality and cluster results are most of the time strongly affected by the chosen method and its performance is very dependent on the underlying probabilistic model which the method implicitly assumes.

For instance, when using $k$-means, we must keep in mind that this method is designed for clustering spherical groups of roughly equal sizes and, thus, it is not reliable for analyzing constellations of groups that depart strongly

Received June 2006; revised June 2006.

[1]Supported in part by Ministerio de Ciencia y Tecnología and FEDER Grant MTM2005-08519-C02-01 and by Consejería de Educación y Cultura de la Junta de Castilla y León Grant PAPIJCL VA102A06.

AMS 2000 subject classifications. Primary 62H3; secondary 62H3.

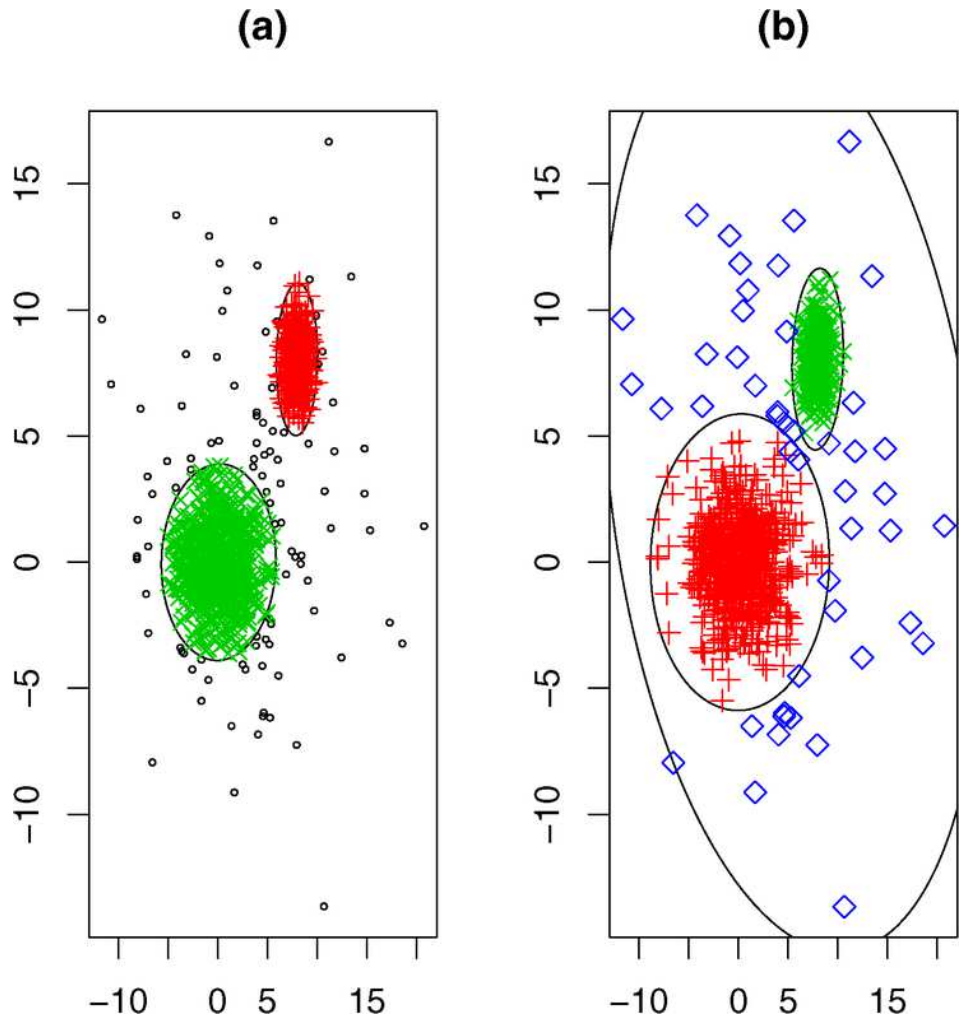

FIG. 1.   (a) *Two groups with 10% of the observation discarded (trimmed points are the small circles).* (b) *Three groups partition with no observations discarded.*

from this assumption. So, in order to understand clustering methods and decide what is more appropriate in a particular case, it is interesting to construct feasible models and develop suitably tailored methods for them.

Determining appropriate models for clustering is even more important when noisy data or outliers are present. Without specifying a model, what we understand by an observation following an "anomalous" behavior is not clear. For instance, it is difficult to decide when a set of very scattered observations should be considered as an extra proper group or merely as a background noise to be discarded (see Figure 1). Additionally, it is not obvi-ous if a small group of tightly joined outliers should be considered as a proper group instead of a contamination phenomenon. Finally, note that the precise detection of the outliers is an important task due to the serious troubles they introduce in standard clustering procedures (see, e.g., García-Escudero and Gordaliza [12] and Hennig [19]) as well as the appealing interest that outliers could have by themselves after explaining why they depart from general behavior.

Two general model-based approaches which provide a theoretically well-based clustering criterion in presence of outliers are (see Bock [2]) the *mixture modeling* and the *trimming approach*. To the first category belongs, say, the work by Fraley and Raftery [8], that considers mixture fittings with the addition of a mixture component accounting for the "noise," or McLachlan and Peel [23] that resorts to mixtures of $t$ distributions. In this paper we are concerned with the trimming approach, previously introduced in Cuesta-Albertos, Gordaliza and Matrán [4] and followed by recent proposals by Gallegos [9, 10] and Gallegos and Ritter [11] (see also García-Escudero, Gordaliza and Matrán [14] and [15]). Notice that a "crisp" 0–1 approach is usually adopted in trimming approaches while some groups' ownership probabilities are generally returned by mixture modeling. Also, while mixture modeling tries to fit the outlying observations in the model, the trimming approach attempts to discard them completely. The methodology presented in this paper falls within the category of trimming approach methods and all the comparisons will be made within this category.

To know how to perform the trimming in cluster analysis is not straightforward because there exist no privileged directions for searching outlying values and, most of the time, we even need to remove observations which fall between the groups ("bridge" data points). The first attempt of trimming in clustering, through an "impartial" approach, appeared in [4] as a modification of the $k$-means method. Moreover, [12] shows that the impartial trimming provides better results in terms of robustness than the consideration of different penalty functions in the $k$-means method (e.g., $k$-medoids).

The use of trimmed $k$-means involves a considerable drawback because it implicitly assumes the same spherical covariance matrix for the groups (as classical $k$-means does). The extension in [11] through the trimmed determinant criterion allows for a general expression of a common covariance matrix. Moreover, [11] also introduces there a statistical clustering model with outliers called the *spurious-outlier model* extending the usual statistical clustering setup (Mardia, Kent and Bibby [20]) to include the presence of a proportion $\alpha$ of noise. This point of view leads to the consideration of the clustering method via maximum likelihood that we pursue in this paper.

Unfortunately, the heterogeneous robust clustering problem (where different groups' covariance matrices are admitted) is notably harder. The proposed objective function is now unbounded and the different "scales"

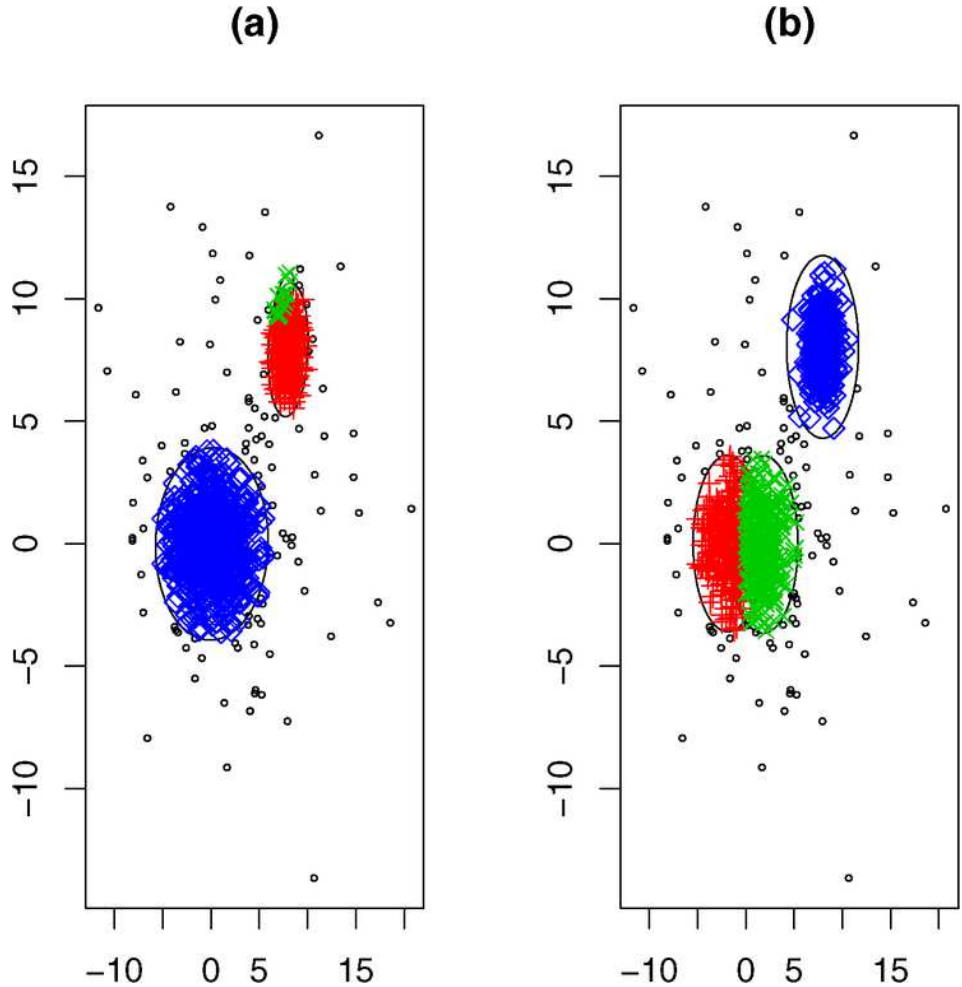

Fig. 2.   *An unrestricted solution for the same data set in Figure 1 appears in* (a) *when*
*k = 3 and α = 0.1. Compare with a restricted solution, also when k = 3 and α = 0.1, in*
(b).

complicates the global ordering of the observations around their closest
centers through Mahalanobis distances (see García-Escudero and Gordaliza
[13]). This motivates that unrestricted algorithms often find small clusters
of points either grouped or almost lying in a lower dimensional space [Figure 2(a)]. Adding some kind of restriction could allow us to obtain more
informative partitions [Figure 2(b)].

A possible way of adding restriction has been considered in Gallegos [9, 10]
by normalizing the covariances to have unit determinant in the steps followed
in the algorithms there. This idea works nicely when the groups have similar scales, but it does not work so well when very different scales are involved. This normalization can be too restrictive and it seems more adequate to incorporate the restrictions directly in the problem statement instead of (artificially) appearing in the algorithm.

To address these difficulties, we introduce in our proposal constraints on the covariance matrices eigenvalues-ratio. A constant $c$ will control the strength of the restrictions allowing a wide range of clustering problems.

Another difficulty arises under the presence of different sizes for the underlying groups. Our proposal also includes, in a successful way, the consideration of different groups' weights to handle this difficulty.

Existence results for both, the sample and the population problem, as well as the consistency of the sample maximizers to the population ones under mild assumptions are shown in Section 2. The proofs are sketched in the Appendix stressing on the importance of the eigenvalue restrictions to achieve these results.

In Section 3, we propose a feasible algorithm (TCLUST) for approximately solving the sample version of the problem. It may be seen as a classification EM-algorithm (Celeux and Govaert [3]) where a kind of "concentration" step as in fast-MCD algorithm (Rousseeuw and van Driessen [25]) is also applied. The eigenvalues-ratio restrictions will be imposed by solving a restricted least squares problem. Dykstra's algorithm in [6] may be applied for addressing that problem. Finally, in Section 4, we include a simulation study showing the gain provided by the proposed method with respect to other trimming proposals.

**2. Robust clustering and eigenvalues-ratio restrictions.** We will consider throughout the paper a data set $\{x_1, \ldots, x_n\}$ in the Euclidean space $\mathbb{R}^p$. By $f(\cdot; \mu, \Sigma)$, we will denote the probability density function (p.d.f.) of the $p$-variate normal distribution with mean $\mu$ and covariance matrix $\Sigma$.

Under the spurious-outlier model, introduced in [11], the likelihood function is given by

$$(2.1) \qquad \left[ \prod_{j=1}^{k} \prod_{i \in R_j} f(x_i; \mu_j, \Sigma) \right] \left[ \prod_{i \notin R} g_i(x_i) \right]$$

with $R = \bigcup_{j=1}^{k} R_j$ and $\#R = n - [n\alpha]$, and where the parameter $k$ denotes the total number of groups, $R_j$ contains the indexes of the "regular" observations assigned to group $j$ and the remaining observations are considered spurious and obtained from some $g_i$'s, p.d.f.s in $\mathbb{R}^p$.

If $\Sigma = \sigma^2 I$ is chosen in (2.1), then we would be performing the trimmed $k$-means method. An algorithm in the spirit of the fast-MCD (both coincide when $k = 1$) is provided in [11] for approximately maximizing (2.1).

Our modification of the "spurious-outlier" model considers different scatter matrices $\Sigma_j$'s and assumes the presence of some underlying group weights $\pi_j$'s, with $\sum_{j=1}^{k} \pi_j = 1$. This leads to the maximization of

$$(2.2) \qquad \left[ \prod_{j=1}^{k} \prod_{i \in R_j} \pi_j f(x_i; \mu_j, \Sigma_j) \right] \left[ \prod_{i \notin R} g_i(x_i) \right],$$

with $R = \bigcup_{j=1}^{k} R_j$ and $\#R = n - [n\alpha]$. Additionally, restrictions on the eigenvalues of the $\Sigma_j$'s matrices will be later introduced.

As in [11], we can avoid the nonregular contribution to the previous maximization problem when the $g_i$'s satisfy the condition

$$(2.3) \quad \arg\max_{\mathcal{R}} \max_{\mu_j, \Sigma_j} \prod_{j=1}^{k} \prod_{i \in R_j} \pi_j f(x_i; \mu_j, \Sigma_j) \subseteq \arg\max_{\mathcal{R}} \prod_{i \notin \bigcup_{j=1}^{k} R_j} g_i(x_i),$$

where $\mathcal{R}$ stands for the set of all partitions of the indexes $\{1, \dots, n\}$ onto $k$ groups of regular observations, $R$, and a group containing the nonregular ones, with $\#R = n - [n\alpha]$. Note that the right-hand side in condition (2.3) only involves the nonregular observations and does not depend on the partition of the regular ones. Therefore, it simply means that any set of nonregular observations in every optimal partition maximizing (2.2) could be also obtained as a subset of $[n\alpha]$ elements of the sample maximizing the likelihood corresponding to the noise. This condition easily holds under reasonable assumptions for the $g_i$'s whenever the nonregular observations may be seen as merely "noise." For instance, examples for $g_i$'s shown in [9, 10] and [11] can be trivially considered here. We refer the interested reader to these papers for details.

A better statement of our problem is obtained by introducing *assignment functions* $z_j$, $j = 0, 1, \dots, k$. For every point $x$ in $\mathbb{R}^p$ (not only the sample observations $x_i$'s are classified), let us define $z_j(x) = 1$ whenever $x$ is assigned to the class $R_j$, $j = 1, \dots, k$, or $z_0(x) = 1$ if it is being trimmed off. Through these functions, assuming that the $g_i$'s may be omitted, we can raise again the problem in (2.2) to the maximization of

$$\prod_{i=1}^{n} \left[ \prod_{j=1}^{k} \pi_j^{z_j(x_i)} f(x_i; \mu_j, \Sigma_j)^{z_j(x_i)} \right],$$

where $z_j$ are 0–1 functions defined in the whole sample space verifying $\sum_{j=0}^{k} z_j(x_i) = 1$ and $\sum_{i=1}^{n} z_0(x_i) = [n\alpha]$. This statement of the problem, taking logarithms, leads to the following general one.

**Robust clustering problem:** Given a probability measure $P$, maximize

$$(2.4) \qquad E_P \left[ \sum_{j=1}^{k} z_j(\cdot)(\log \pi_j + \log f(\cdot; \mu_j, \Sigma_j)) \right],$$

in terms of the assignment functions

$$z_j : \mathbb{R}^p \mapsto \{0, 1\} \qquad \text{such that } \sum_{j=0}^{k} z_j = 1 \text{ and } E_P z_0(\cdot) = \alpha,$$

and the parameters $\theta = (\pi_1, \ldots, \pi_k, \mu_1, \ldots, \mu_k, \Sigma_1, \ldots, \Sigma_k)$ corresponding to weights $\pi_j \in [0, 1]$, with $\sum_{j=1}^{k} \pi_j = 1$, mean vectors $\mu_j \in \mathbb{R}^p$ and symmetric positively definite $p \times p$-matrices $\Sigma_j$, $j = 1, \ldots, k$.

If $P_n$ stands for the empirical measure, $P_n = 1/n \sum_{i=1}^{n} \delta_{\{x_i\}}$, by replacing $P$ by $P_n$, we recover the original sample problem [notice that, perhaps, $E_{P_n} z_0(\cdot) = \alpha$ cannot be exactly achieved but this familiar fact will not be important in our reasonings].

Our restrictions on the eigenvalues of the covariance matrices may be seen as an extension of those introduced by Hathaway [18] for univariate data. They avoid the singularities introduced by the possibility of very different $\Sigma_j$'s.

**(ER) Eigenvalues-ratio restrictions:** We fix a constant $c \geq 1$ such that

$$M_n / m_n \leq c$$

for

$$M_n = \max_{j=1,\ldots,k} \max_{l=1,\ldots,p} \lambda_l(\Sigma_j) \quad \text{and} \quad m_n = \min_{j=1,\ldots,k} \min_{l=1,\ldots,p} \lambda_l(\Sigma_j),$$

$\lambda_l(\Sigma_j)$ being the eigenvalues of the matrices $\Sigma_j$, $l = 1, \ldots, p$ and $j = 1, \ldots, k$. The set of $\theta$'s which obey this condition is denoted by $\Theta_c$.

Note that $c = 1$ produces the strongest possible restriction. In this case, the proposed method may be viewed as a trimmed $k$-means method with weights. However, the main advantage of this approach relies on the fact that the parameter $c$ allows us to achieve certain (controlled) freedom in how we want to handle the different scattering of the groups.

Figures 1 and 2 show the results of the application of the proposed methodology (by using the TCLUST algorithm described in Section 3) to a data set made up of 3 bivariate Gaussian clusters where the most scattered one accounts for 10% of the data. The result when $k = 2$, $\alpha = 0.1$ and $c = 5$ appears in Figure 1(a). The result there is not very dependent on $c$ as long as the two (main) groups are not too different in their eigenvalues once the most scattered group was trimmed off. The values $k = 3$ and $\alpha = 0$ are considered in Figure 1(b), with a large value for $c$ ($c = 50$) which allows for the presence of the more scattered group. The values $k = 3$ and $\alpha = 0.1$ were applied in Figure 2. A rather large $c$ (unrestricted problem) was chosen in Figure 2(a) while a small $c = 1$ (restricted problem) was considered in Figure 2(b).

To exclude in the subsequent analysis those probability distributions obviously unappropriate for the introduced approach, we will assume on the underlying distribution $P$ the following mild condition. (It trivially holds for absolutely continuous distribution or for empirical measures corresponding to a sample large enough from an absolutely continuous distribution.)

**(PR)** The distribution $P$ is not concentrated on $k$ points after removing a probability mass equal to $\alpha$.

To conclude this section, we will notably simplify our problem through an adequate reformulation that leads to expressing the assignment functions $z_j$'s only in terms of $\theta$. This will also be a keystone for deriving our algorithm to solve the sample counterpart of the problem.

Given $\theta \in \Theta_c$, we consider *discriminant functions* defined as

$$
\begin{aligned}
(2.5) \qquad & D_j(x;\theta) = \pi_j f(x;\mu_j,\Sigma_j) \quad \text{and} \\
& D(x;\theta) = \max\{D_1(x;\theta),\ldots,D_k(x;\theta)\}.
\end{aligned}
$$

Note that these are familiar functions in the application of Bayes' rules in discriminant analysis. These functions will also serve to provide an "outlyingness" measure of the observations.

Using previous definitions, for a given $\theta$ and a probability measure $P$, we consider the distribution function of $D(\cdot;\theta)$ and its corresponding $\alpha$-quantile:

$$
(2.6) \quad G(u;\theta,P) := P(D(\cdot;\theta) \leq u) \quad \text{and} \quad R(\theta,P) := \inf_u\{G(u;\theta,P) \geq \alpha\}.
$$

With this notation, we have the following straightforward characterization for the assignment functions:

PROPOSITION 1. *The robust clustering problem can be simplified, using the discriminant functions (2.5), to the maximization in $\theta$ of*

$$
(2.7) \qquad \theta \mapsto L(\theta,P) := E_P\left[\sum_{j=1}^k z_j(\cdot;\theta)\log D_j(\cdot,\theta)\right],
$$

*where the assignment functions are obtained from $\theta$ as*

$$
z_j(x;\theta) = I\{x : \{D(x;\theta) = D_j(x;\theta)\} \cap \{D_j(x;\theta) \geq R(\theta,P)\}\}
$$

*and*

$$
z_0(x;\theta) = 1 - \sum_{j=1}^k z_j(x;\theta).
$$

That is, we assign $x$ to the class $j$ with the largest discriminant function value $D_j(x;\theta)$ or $x$ is trimmed off when all the $D_j(x;\theta)$'s [and consequently, $D(x;\theta)$] are smaller than $R(\theta, P)$. (In order to break ties in the discriminant function values, the lexicographical ordering could be applied.)

The relevant mathematical results to be considered are given in the following propositions.

PROPOSITION 2 (Existence). *If* (PR) *holds for distribution* $P$, *then there exists some* $\theta \in \Theta_c$ *such that the maximum of* (2.4) *under* (ER) *is achieved.*

By examining the proof of previous result, we can see that although we have admitted weights $\pi_j = 0$, this is not a drawback when taking $\log \pi_j$ because in this case $z_j(\cdot; \theta) \equiv 0$ and then the set $\{x : z_j(x; \theta) = 1\}$ is empty. The presence of groups with zero weight does actually happen in practice. For instance, when $k = 2$, $c = 1$, $\alpha = 0$ and $P$ is the $N(0, 1)$ distribution in the real line, we can see that $\theta = (\pi_1, \pi_2, \mu_1, \mu_2, \sigma_1^2, \sigma_1^2) = (1, 0, 0, \mu_2, 1, 1)$ is the optimal solution for every $\mu_2 \in \mathbb{R}$ and $c \geq 1$.

PROPOSITION 3 (Consistency). *Assume that* $P$ *has a strictly positive density function and that* $\theta_0$ *is the unique maximum of* (2.4) *under* (ER). *If* $\theta_n \in \Theta_c$ *denotes a sample version estimator based on the empirical measure* $P_n$, *then* $\theta_n \to \theta_0$ *almost surely.*

REMARK 1. Notice that a uniqueness condition is needed in order to establish the consistency result. Unfortunately, this property does not always hold. For instance, think of a symmetric mixture $P$ in the real line with two well-separated modes, a high trimming level and $k = 1$. That uniqueness property was already needed for establishing the same consistency result for the trimmed $k$-means and, even in this simpler case, the statement of general uniqueness results was difficult (see Remark 4.1 in [14]). However, as in the trimmed $k$-means problem, we believe that it is quite rare to find distributions where the uniqueness fails when dealing with "reasonable" data for clustering and when parameters $k$ and $\alpha$ have been properly chosen.

**3. The TCLUST algorithm.** The empirical problem presented in Section 2 has a very high computational complexity. An exact algorithm seems to be not feasible even for moderate sample sizes. Thus, the existence of adequate algorithms for approximately solving the sample problem is as important as the procedure itself. With this in mind, we propose the TCLUST algorithm (an R-code implementation is available at http://www.eio.uva.es/~langel/software), an EM-principle based algorithm, intended to search for approximate solutions. The EM algorithm (Dempster, Laird and Rubin [5]) is the usual method for obtaining a solution to the mixture likelihood problem. Here, as we follow a "crisp" approach where each point is uniquely assigned to one cluster, a classification EM approach (Celeux and Govaert [3]) is preferable. Moreover, as trimmed observations are allowed, the rationale behind the fast-MCD [25] and behind the trimmed $k$-means algorithm [15] will also underly.

The TCLUST algorithm may be described as follows:

1. Randomly select starting values for the centers $m_j^0$'s, the covariance matrices $S_j^0$'s and the weights of the groups $p_j^0$'s for $j = 1, \ldots, k$.
2. From the $\theta^l = (p_1^l, \ldots, p_k^l, m_1^l, \ldots, m_k^l, S_1^l, \ldots, S_k^l)$ returned by the previous iteration:
   2.1. Obtain $d_i = D(x_i, \theta^l)$ for the observations $\{x_1, \ldots, x_n\}$ and keep the set $H$ having the $[n(1 - \alpha)]$ observations with largest $d_i$'s.
   2.2. Split $H$ into $H = \{H_1, \ldots, H_k\}$ with $H_j = \{x_i \in H : D_j(x_i, \theta^l) = D(x_i, \theta^l)\}$.
   2.3. Obtain the number of data points $n_j$ in $H_j$ and their sample mean and sample covariance matrix, $m_j$ and $S_j$, $j = 1, \ldots, k$.
   2.4. Consider the singular-value decomposition of $S_j = U_j' D_j U_j$ where $U_j$ is an orthogonal matrix and $D_j = \text{diag}(\Lambda_j)$ is a diagonal matrix (with diagonal elements given by the vector $\Lambda_j$). If the full vector of eigenvalues $\Lambda = (\Lambda_1, \ldots, \Lambda_k)$ does not satisfy the eigenvalues-ratio restriction, obtain (for instance) through Dykstra's algorithm a new vector $\tilde{\Lambda} = (\tilde{\Lambda}_1, \ldots, \tilde{\Lambda}_k)$ obeying the (ER) restriction and with $\|\tilde{\Lambda} - \Lambda^{-1}\|^2$ being as smaller as possible. ($\Lambda^{-1}$ denotes the vector made up by the inverse of the elements of the vector $\Lambda$.) Notice that the (ER) restriction for $\Lambda$ corresponds exactly to the same (ER) restriction applied to $\Lambda^{-1}$.
   2.5. Update $\theta^{l+1}$ by using:
   - $p_j^{l+1} \leftarrow n_j/[n(1 - \alpha)]$,
   - $m_j^{l+1} \leftarrow m_j$,
   - $S_j^{l+1} \leftarrow U_j' \tilde{D}_j U_j$ and $\tilde{D}_j = \text{diag}(\tilde{\Lambda}_j)^{-1}$.
3. Perform $F$ iterations of the process described in step 2 (moderate values for $F$ are usually enough) and compute the evaluation function $L(\theta^F; P_n)$ using expression (2.7).
4. Start from step 1 several times, keeping the solutions leading to minimal values of $L(\theta^F, P_n)$ and fully iterate them to choose the best one.

The computed (E-step) "a posteriori" probabilities, $D_j(x_i, \theta^l) = p_j f(x_i; m_j, S_j)$, are converted to a discrete classification leaving unassigned the proportion $\alpha$ of observations which are the hardest to classify. It is easy to see that this leads us to an optimal assignment. We later obtain a new $\theta^{l+1}$ by maximizing (M-step) the conditional expectation. Proposition 4 guarantees that the presented algorithm can be applied for performing this maximization. Notice that the obtention of the optimal scatter matrices is decomposed into the search of the corresponding optimal eigenvalues and eigenvectors. For every choice of eigenvalues, the best eigenvectors choice simply follows from those derived from the sample covariance matrices of the observations in each group. This decomposition is somehow similar to that considered in Gallegos' proposal, where "shapes" and "scales" were separately handled.

Seeing $D(x_i, \theta^l)$ as an inverse outlyingness measure for the observation $x_i$ with respect to the choice of $\theta^l$, then step 2 may be seen as a concentration step. [13] analyzes some other attempts for extending the concentration step principle to the heterogeneous robust clustering setup.

Recall that the random initialization scheme (step 1) and the final refinement (step 4) will be very important as happened in the fast-MCD algorithm or in Maronna [21]. For initializing the procedure in step 1, we have seen that simply randomly choosing $k$ sample data points for the centers, $k$ identity matrices for the covariances and the same weights for the groups (equal to $1/k$) provide reasonably starting values in most of the cases.

With respect to the eigenvalues-ratio restriction, we would need $\Lambda = (\Lambda_1, \ldots, \Lambda_k)$ with $\Lambda_j = (\lambda_{1,j}, \ldots, \lambda_{p,j})$ belonging to the cone

$$(3.1) \quad \mathcal{C} = \{(\Lambda_1, \ldots, \Lambda_k) \in \mathbb{R}^{p \times k} : \lambda_{u,v} - c \cdot \lambda_{r,s} \le 0 \text{ for all } (u,v) \ne (r,s)\}.$$

If $\Lambda \notin \mathcal{C}$, we must replace $\Lambda^{-1}$ by $\tilde{\Lambda} \in \mathcal{C}$ with minimal $\|\tilde{\Lambda} - \Lambda^{-1}\|^2$. Dykstra's algorithm serves to approximately solve that restricted least squares problem as long as $\mathcal{C}$ is the intersection of the several closed convex cones

$$\mathcal{C}_h = \{(\Lambda_1, \ldots, \Lambda_k) \in \mathbb{R}^{p \times k} : \lambda_{u,v} - c \cdot \lambda_{r,s} \le 0\} \qquad \text{for } h = (u, v, r, s),$$

by resorting to iterative projections onto the individual cones $\mathcal{C}_h$'s. (Notice that the projections onto the cones $\mathcal{C}_h$ are very fast to obtain.) Thus, a fixed number of individual projections may be done retaining the best attained solution after these iterations and satisfying the restrictions. Alternatively, quadratic programming based solutions (see, e.g., Goldfarb and Idnani [17]) for that constrained minimization may be explored.

The next result formalizes the appropriateness of the TCLUST algorithm.

PROPOSITION 4. *If the sets $H_j = \{x_i : z_j(x_i) = 1\}$, $j = 1, \ldots, k$, are kept fixed, the maximum of (2.4) for $P = P_n$ can be obtained through the following steps:*

(i) *Fixed $\mu_j$ and $\Sigma_j$, the best choice of $\pi_j$ is $\pi_j = n_j/[n(1-\alpha)]$ where $n_j$ is the cardinal of set $H_j$.*

(ii) *Fixed $\Sigma_j$ and the optimal values for $\pi_j$ given in* (i), *the best choice for $\mu_j$ is the sample mean $m_j$ of the observations in $H_j$.*

(iii) *Fixed the eigenvalues for the matrix $\Sigma_j$ and the optimum values given in* (i) *and* (ii) *for $\pi_j$ and $\mu_j$, the best choice for the set of unitary eigenvectors are the unitary eigenvectors of the sample covariance matrix $S_j$ of the observations in $H_j$.*

(iv) *With the optimal selections made in* (i), (ii) *and* (iii), *the best choice for the eigenvalues corresponds to the inverse of the projection of the vector containing the inverse of the eigenvalues onto the cone $\mathcal{C}$.*

PROOF. Once the $z_j(x_i)$ for $i = 1, \ldots, n$ and $j = 0, \ldots, k$ are known values, the expression (2.4) can be written as

$$(3.2) \qquad \sum_{j=1}^{k} \left[ n_j \log \pi_j + \sum_{x_i \in H_j} \log f(x_i; \mu_j, \Sigma_j) \right],$$

and the assertions (i) and (ii) trivially hold.

Considering these optimal values for $\pi_j$ and $\mu_j$, together with the cyclic property of the trace, the maximization of (3.2) simplifies to the minimization of

$$\sum_{j=1}^{k} [\log |\Sigma_j| + \operatorname{trace}(\Sigma_j^{-1} S_j)].$$

The matrices $S_j$ and $\Sigma_j$ can be decomposed into $S_j = U_j' D_j U_j$ and $\Sigma_j = V_j' E_j V_j$, where $D_j = \operatorname{diag}(\Lambda_j)$ and $E_j = \operatorname{diag}(\Xi_j)$ are diagonal matrices with $\Lambda_j = (\lambda_{1,j}, \ldots, \lambda_{p,j})$ and $\Xi_j = (\xi_{1,j}, \ldots, \xi_{p,j})$, and $U_j$ and $V_j$ are orthogonal matrices. So, as $\log |\Sigma_j| = \log |E_j|$ and the eigenvalues $E_j$ were fixed, the previous minimization problem can be further simplified to that of

$$\sum_{j=1}^{k} \operatorname{trace}(\Sigma_j^{-1} S_j) = \sum_{j=1}^{k} \operatorname{trace}(E_j^{-1}(U_j V_j')' D_j (U_j V_j'))$$

(the cyclic property of the trace is again applied). Denote $T_j = U_j V_j'$ and rewrite

$$(3.3) \qquad \operatorname{trace}(E_j^{-1} T_j' D_j T_j) = \sum_u \sum_v \frac{\lambda_{u,j}}{\xi_{v,j}} \cdot t_{uv,j}^2,$$

where $t_{uv,j}$ denotes the element $(u, v)$ of the matrix $T_j$. $T_j$ is an orthogonal matrix, so we have that $\sum_u t_{uv,j}^2 = 1$ and $\sum_v t_{uv,j}^2 = 1$. Therefore, the minimization of (3.3) may be seen as a linear programming problem like

$$\min \sum_{u,v} c_{u,v} \cdot x_{u,v} \quad \text{subject to} \quad \sum_u x_{u,v} = 1, \sum_v x_{u,v} = 1 \text{ and } x_{u,v} \geq 0,$$

with known coefficients $c_{u,v}$ [notice that $\lambda_{u,j}/\xi_{u,j}$ are fixed coefficients because $\lambda_{u,j}$ depends on the data set at hand and the $\xi_{u,j}$ are supposed known quantities in (iii)]. Although fractional solutions are possible, these solutions will never be basic feasible ones due to the particular statement of the linear programming problem (see, e.g., Papadimitriou and Steiglitz [24], page 249). Consequently, the optimal solution corresponds to a "real matching" where the optimal $t_{u,v}^2$ are 0 or 1. Thus, $T_j$ is a permutation matrix product of the orthogonal matrices $U_j$ and $V_j'$. It is quite easy to see that the columns of the matrices $U_j$ and $V_j$ must provide the same set of unitary eigenvectors and, thus, the assertion (iii) is proven.

By applying (i), (ii) and (iii), we finally need to search for a vector $\Xi = (\Xi_1, \ldots, \Xi_k)$ and $\Xi_j = (\xi_{1,j}, \ldots, \xi_{p,j})$ minimizing

$$(3.4) \qquad \sum_{j=1}^{k} \sum_{i=1}^{p} \left( \log \xi_{i,j} + \frac{\lambda_{i,j}}{\xi_{i,j}} \right) = \sum_{j=1}^{k} \sum_{i=1}^{p} (-\log \tilde{\lambda}_{i,j} + \lambda_{i,j} \cdot \tilde{\lambda}_{i,j}),$$

with $\tilde{\lambda}_{i,j} = 1/\xi_{i,j}$. As (3.4) is a convex function on the $\tilde{\lambda}_{i,j}$ and its unrestricted minimum is attained when $\tilde{\lambda}_{i,j} = \lambda_{i,j}^{-1}$, the minimization of (3.4) under the eigenvalues-ratio restriction posed by (3.1) leads us to the optimal choice of $\tilde{\Lambda}$ with minimal $\|\tilde{\Lambda} - \Lambda^{-1}\|^2$ and $\tilde{\Lambda} \in \mathcal{C}$.   $\square$

REMARK 2.   Alternative methods can be defined by imposing restrictions on the ratio between covariance determinants instead of controlling eigenvalues. Gallegos [9] and [10] proposal scales the covariance matrices to have determinant ratio equal to 1 in the algorithm. Maronna and Jacovkis [22], in the untrimmed case $\alpha = 0$, consider that normalization as the only reliable "distance" for clustering multivariate data. Here, the proposed algorithm can be easily adapted for handling restrictions of this type. In this case, the cone would be

$$\mathcal{C}' = \{(\sigma_1, \ldots, \sigma_k) \in \mathbb{R}^k : \sigma_u - c \cdot \sigma_v \le 0 \text{ for all } u \ne v\},$$

where the factorization in step 2.4 of the previous algorithm is $S_j = \sigma_j \cdot U_j$ with $|U_j| = 1$ and $\sigma_j = |S_j|^{1/p}$. If $c = 1$ in $\mathcal{C}'$, we would obtain an analogous to Gallegos' proposal with group weights.

Other procedures which have been used for avoiding pathological solutions in the heterogeneous robust clustering problem are based on adding different types of parameterizations for the covariance matrices (see, e.g., Scott and Symons [26] or Banfield and Raftery [1]). Although that possibility has not been considered here, we believe that similar ideas (based on relaxing those parameterizations) could be interesting.

REMARK 3.   The proper determination of parameters $\alpha$, $k$ and $c$ is not an easy problem in general. Users of cluster analysis methods sometimes have initial guesses of suitable values for these parameters, but many times these are completely unknown. The careful analysis of the objective functions when moving $k$ and $\alpha$ provided useful information for choosing $k$ and $\alpha$ in [15]. The objective function for the trimmed $k$-means method always improves when increasing $k$ (see Lemma 2.2 in [4]). Here, the possible existence of groups with $\pi_j = 0$ would imply that the value of the objective function does not necessarily improve when increasing $k$. However, this property could even be more interesting in order to develop techniques for choosing $k$ because a $\pi_j$ close to 0 suggests that an smaller $k$ could be needed.

Moreover, as an anonymous referee suggested to us, we can make use of Bayes factors as in Van Aelst et al. [27] in order to know how well the observation $x_i$ is integrated in the cluster in which it was assigned. If $D_{(1)}(x_i; \theta) \leq \cdots \leq D_{(k)}(x_i; \theta)$, the Bayes factor for a nontrimmed observation $x_i$ is defined as $\mathrm{BF}(i) = \log(D_{(k-1)}(x_i; \theta) / D_{(k)}(x_i; \theta))$. Notice that the smaller the Bayes factor is the better is the assignment to its corresponding cluster. The existence of clusters with many observation with large Bayes factors (close to 0) suggests that perhaps an improper choice for $c$ was made. Additionally, we can introduce Bayes factors for the trimmed observations as $\mathrm{BF}(i) = \log(D_{(k)}(x_i, \theta) / R(\theta, P_n))$ to measure the strength of the consideration of the trimmed data point $x_i$ as an outlier.

## 4. A simulation study.

A simulation study has been carried out to compare the performance of the proposed robust clustering method with respect to other trimming approaches in the literature. Several data sets of size $n = 2000$ have been generated. Each data set consists of three simulated $p$-dimensional normally distributed clusters with centers $\mu_1 = (0, 8, 0, \ldots, 0)'$, $\mu_2 = (8, 0, 0, \ldots, 0)'$ and $\mu_3 = (-8, -8, 0, \ldots, 0)'$ and covariance matrices

$$\Sigma_1 = \mathrm{diag}(1, a, 1, \ldots, 1),$$
$$\Sigma_2 = \mathrm{diag}(b, c, 1, \ldots, 1)$$

and

$$\Sigma_3 = \begin{pmatrix} \begin{array}{cc} d & e \\ e & f \end{array} & \mathbf{0} \\ \hline \mathbf{0} & I \end{pmatrix}.$$

The constants $a, b, c, d, e$ and $f$ serve to control the true differences between the eigenvalues of the groups' covariance matrices. This leads us to the consideration of the following cases:

(M1) $(a, b, c, d, e, f) = (1, 1, 1, 1, 0, 1)$: *Spherical equally scattered groups.*

(M2) $(a, b, c, d, e, f) = (5, 1, 5, 1, 0, 5)$: *Not spherical but the same covariance matrices for the groups.*

(M3) $(a, b, c, d, e, f) = (5, 5, 1, 3, -2, 3)$: *Different covariance matrices but the same scale (equal determinant).*

(M4) $(a, b, c, d, e, f) = (1, 20, 5, 15, -10, 15)$: *Groups with different scales.*

(M5) $(a, b, c, d, e, f) = (1, 45, 30, 15, -10, 15)$: *Groups with different scales and two of them with a severe overlap.*

We consider 1800 "regular" data points with two different group proportions. We also generate uniformly distributed data points in a parallelogram defined by the coordinatewise ranges of the regular data points. Using an acceptance–rejection algorithm, only points having squared Mahalanobis distances from $\mu_1$, $\mu_2$ and $\mu_3$ (using $\Sigma_1$, $\Sigma_2$ and $\Sigma_3$) greater than $\chi^2_{p,0.975}$ are finally considered until reaching an amount of 200 outliers.

The following approaches searching for $k = 3$ groups and with a trimming proportion $\alpha = 0.1$ are tried:

(TkM) *Trimmed k-means (specially aimed to the case* M1).

(G&R) *Gallegos and Ritter's method (specially aimed to the case* M2).

(G) *Gallegos' proposal (specially aimed to the case* M3).

(TCLUST) *The presented algorithm with an eigenvalues-ratio restriction when* $c = 50$.

The same number of random initializations and concentration steps are taken for all the methods.

Table 1 shows the average proportion of misclassified observations for $B = 100$ independent random samples of size 2000 when $p = 2$ and 6 and the group weights satisfy proportions 1:1:1 ("equal") and 1:2:2 ("unequal"). The numbers within parenthesis in the Table 1 show the proportion of outliers wrongly determined as nonoutliers and vice versa.

Notice that all the methods work nicely under the underlying model in which they are specially aimed. However, the proposed eigenvalues-ratio restriction method is the only method which is able to cope with the mixtures with very different scales (mixtures M4 and M5), and it seems to be less affected in the unequal groups' size case. Figure 3 shows the result of these four analyzed procedures applied to the same data set generated by the simulation scheme M5 when $p = 2$ and unequal weights. The TCLUST seems to be the only one that is able to distinguish between the least and the most scattered groups even in this rather overlapped case.

## APPENDIX: PROOFS OF EXISTENCE AND CONSISTENCY

**A.1. Existence.** The existence of solutions for our problem can be obtained through a standard argument starting by considering a sequence $\{\theta_n\}_{n=1}^{\infty} = \{(\pi_1^n, \ldots, \pi_k^n, \mu_1^n, \ldots, \mu_k^n, \Sigma_1^n, \ldots, \Sigma_k^n)\}_{n=1}^{\infty}$ such that

$$(A.1) \qquad \lim_{n \to \infty} L(\theta_n, P) = \sup_{\theta \in \Theta_c} L(\theta, P) = M > -\infty$$

TABLE 1
*Misclassification rates for the simulation study*

| Weights | Dimen. | Model | TkM | G&R | G | TCLUST |
|---------|--------|-------|-----|-----|---|--------|
| Equal | $p = 2$ | M1 | 0.011 (0.011) | 0.011 (0.011) | 0.012 (0.012) | 0.012 (0.012) |
| | | M2 | 0.037 (0.037) | 0.016 (0.016) | 0.016 (0.016) | 0.016 (0.016) |
| | | M3 | 0.036 (0.036) | 0.035 (0.035) | 0.015 (0.014) | 0.015 (0.014) |
| | | M4 | 0.079 (0.057) | 0.059 (0.049) | 0.042 (0.031) | 0.020 (0.019) |
| | | M5 | 0.129 (0.046) | 0.131 (0.046) | 0.122 (0.042) | 0.043 (0.022) |
| | $p = 6$ | M1 | 0.008 (0.008) | 0.008 (0.008) | 0.008 (0.008) | 0.009 (0.009) |
| | | M2 | 0.035 (0.035) | 0.012 (0.012) | 0.012 (0.012) | 0.012 (0.012) |
| | | M3 | 0.032 (0.032) | 0.018 (0.017) | 0.011 (0.011) | 0.011 (0.011) |
| | | M4 | 0.094 (0.072) | 0.038 (0.025) | 0.018 (0.016) | 0.015 (0.014) |
| | | M5 | 0.159 (0.077) | 0.119 (0.026) | 0.055 (0.021) | 0.035 (0.015) |
| Unequal | $p = 2$ | M1 | 0.011 (0.011) | 0.011 (0.011) | 0.011 (0.011) | 0.011 (0.011) |
| | | M2 | 0.037 (0.037) | 0.017 (0.017) | 0.017 (0.017) | 0.016 (0.016) |
| | | M3 | 0.036 (0.036) | 0.034 (0.033) | 0.015 (0.015) | 0.014 (0.014) |
| | | M4 | 0.089 (0.059) | 0.069 (0.051) | 0.047 (0.032) | 0.021 (0.019) |
| | | M5 | 0.151 (0.047) | 0.166 (0.048) | 0.147 (0.044) | 0.047 (0.023) |
| | $p = 6$ | M1 | 0.008 (0.008) | 0.009 (0.009) | 0.009 (0.009) | 0.009 (0.009) |
| | | M2 | 0.034 (0.034) | 0.011 (0.011) | 0.012 (0.012) | 0.012 (0.012) |
| | | M3 | 0.033 (0.033) | 0.018 (0.018) | 0.012 (0.011) | 0.011 (0.010) |
| | | M4 | 0.107 (0.074) | 0.053 (0.023) | 0.025 (0.017) | 0.017 (0.015) |
| | | M5 | 0.186 (0.081) | 0.166 (0.024) | 0.078 (0.022) | 0.039 (0.015) |

[the boundedness from below for (A.1) can be easily obtained just considering $\pi_1 = 1$, $\mu_1 = 0$, $\Sigma_1 = I$, and setting the other weights as 0 with arbitrary choices of means and variances].

Since $[0, 1]^k$ is a compact set, we can extract a subsequence from $\{\theta_n\}_n^\infty$ (that will be denoted like the original one) such that

$$(A.2) \qquad \pi_j^n \to \pi_j \in [0, 1] \qquad \text{for } 1 \le j \le k,$$

and satisfying for some $g \in \{0, 1, \ldots, k\}$ (a relabeling could be needed) that

$$(A.3) \qquad \mu_j^n \to \mu_j \in \mathbb{R}^p \qquad \text{for } 0 \le j \le g \quad \text{and} \quad \min_{j > g} \|\mu_j^n\| \to \infty.$$

With respect to the scatter matrices, under (ER), we can also consider a further subsequence verifying one (and only one) of these possibilities:

$$(A.4) \qquad \Sigma_j^n \to \Sigma_j \qquad \text{for } 1 \le j \le k,$$

$$(A.5) \qquad M_n = \max_{j=1,\ldots,k} \max_{l=1,\ldots,p} \lambda_l(\Sigma_j) \to \infty$$

or

$$(A.6) \qquad m_n = \min_{j=1,\ldots,k} \min_{l=1,\ldots,p} \lambda_l(\Sigma_j) \to 0.$$

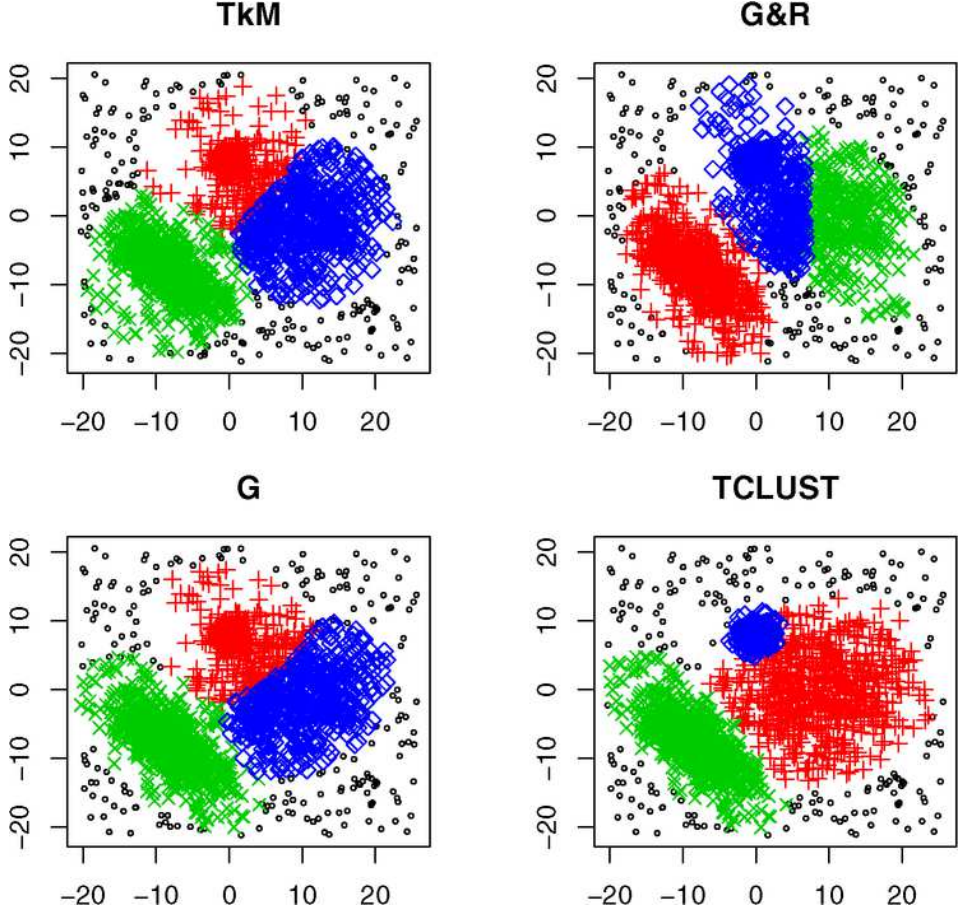

Fig. 3. *Clustering results when $k = 2$ and $\alpha = 0.1$ for a simulated data following the* M5 *scheme in the text with $p = 2$ and unequal weights: Trimmed $k$-means (*TkM*); Gallegos and Ritter (*G&R*); Gallegos (*G*) and the presented algorithm (*TCLUST*) with $c = 50$.*

LEMMA A.1. *If* (ER) *holds and if $P$ satisfies* (PR)*, then only the convergence* (A.4) *is possible.*

PROOF. We will see that (A.5) or (A.6) would imply $\lim_{n \to \infty} L(\theta_n, P) = -\infty$. Let $\lambda_{l,j}^n := \lambda_l(\Sigma_j^n)$ be the eigenvalues, $j = 1, \ldots, k$ and $l = 1, \ldots, p$, of the group covariance matrices and $\|v_{l,j}^n\| = 1$ their associated unitary eigenvectors. Then we have

$$L(\theta_n, P) = E_P\left[\sum_{j=1}^{k} z_j(\cdot; \theta_n)\left(\log \pi_j^n - \frac{p}{2}\log 2\pi - \frac{1}{2}\sum_{l=1}^{p}\log \lambda_{l,j}^n\right.\right.$$

$$-\frac{1}{2}\sum_{l=1}^{p}(\lambda_{l,j}^{n})^{-1}(\cdot-\mu_{j}^{n})'$$

(A.7)
$$\times v_{l,j}^{n}(v_{l,j}^{n})'(\cdot-\mu_{j}^{n})\Big)\Big]$$

$$\leq E_{P}\Bigg[\sum_{j=1}^{k}z_{j}(\cdot;\theta_{n})\Big(\log\pi_{j}^{n}-\frac{p}{2}\log 2\pi$$

$$-\frac{p}{2}\log m_{n}-\frac{1}{2}M_{n}^{-1}\|\cdot-\mu_{j}^{n}\|^{2}\Big)\Bigg].$$

If we assume that $M_{n}\to\infty$ holds, then $m_{n}\to\infty$ by (ER). Thus, we would have that $L(\theta_{n},P)\to-\infty$ leading to a contradiction with (A.1).

Now assume that (A.6) holds. We can guarantee by Lemma A.2 below that if $P$ satisfies (PR), then there exists a constant $h$ such that

(A.8)
$$E_{P}\Bigg[\sum_{j=1}^{k}z_{j}(\cdot;\theta_{n})\|\cdot-\mu_{j}^{n}\|^{2}\Bigg]\geq h>0.$$

Since $\log\pi_{j}^{n}\leq 0$, the fact that $P[z_{1}(\cdot)+\cdots+z_{k}(\cdot)]=1-\alpha$ implies

$$L(\theta_{n},P)\leq(1-\alpha)\Big(-\frac{p}{2}\log 2\pi-\frac{p}{2}\log m_{n}\Big)$$

$$-\frac{1}{2}M_{n}^{-1}E_{P}\Bigg[\sum_{j=1}^{k}z_{j}(\cdot;\theta_{n})\|\cdot-\mu_{j}^{n}\|^{2}\Bigg].$$

Therefore, (ER) and (A.8) give

(A.9)     $$L(\theta_{n},P)\leq(1-\alpha)\Big(-\frac{p}{2}\log 2\pi-\frac{p}{2}\log m_{n}\Big)-\frac{1}{2}(cm_{n})^{-1}h.$$

But this upper-bound in (A.9) tends to $-\infty$ as $m_{n}\to 0$.  □

The following lemma has been applied in the proofs of Lemma A.1.

LEMMA A.2.  *If $P$ satisfies condition* (PR), *then there exists a constant $h>0$ such that inequality* (A.8) *holds.*

PROOF.  The trimmed $k$-means problem was introduced in [4] as the search of $k$ points $m_{1},\ldots,m_{k}$ in $\mathbb{R}^{p}$ and a Borel set $B$ minimizing:

(A.10)
$$\min_{B:\,P(B)\geq 1-\alpha}\min_{m_{1},\ldots,m_{k}}\frac{1}{P(B)}\int_{B}\inf_{1\leq j\leq k}\|x-m_{j}\|\,dP(x).$$

Theorem 3.1 in [4] guarantees the existence of solutions for this problem. Thus, (A.10) attains a minimum value that we denoted by $V_{\alpha,k}$. Now, for every choice of $\theta$, we have

$$E_P\left[\sum_{j=1}^{k} z_j(\cdot;\theta)\|\cdot-\mu_j\|^2\right] \geq E_P\left[\sum_{j=1}^{k} z_j(\cdot;\theta)\inf_{1\leq l\leq k}\|\cdot-\mu_l\|^2\right] \geq (1-\alpha)V_{\alpha,k},$$

because $\bigcup_{j=1}^{k}\{x:z_j(x;\theta)=1\}$ is a Borel set having probability greater or equal than $1-\alpha$. Finally, we can trivially see that $h:=(1-\alpha)V_{\alpha,k}>0$ whenever condition (PR) holds for $P$. $\quad\square$

The next step is to show that whenever the classes in the optimal partition have strictly positive probability masses we can guarantee the convergence of the centers $\mu_j^n$. This result has also key importance in order to understand the role played by the weights $\pi_j$'s in this approach.

Lemma A.3. *When* (ER) *and* (PR) *hold, if every $\pi_j$ in* (A.2) *verifies $\pi_j>0$ for $j=1,\ldots,k$, then $g=k$ in* (A.3).

Proof. If $g=0$, we can take a ball with center 0 and radius big enough $B(0,R)$ such that $P[B(0,R)]>\alpha$. We can thus easily see that

$$E_P\left[\sum_{j=1}^{k} z_j(\cdot;\theta_n)\|\cdot-\mu_j^n\|^2\right]\to\infty,$$

so that $L(\theta_n,P)\to-\infty$ from (A.7). Notice that (ER) is also here applied.

When $g>0$, we prove first that

$$\text{(A.11)}\qquad\qquad E_P\left[\sum_{j=g+1}^{k} z_j(\cdot;\theta_n)\right]\to 0.$$

This arises from the dominated convergence theorem taking into account that the sequence is obviously bounded by $1-\alpha$, and the fact that

$$\text{(A.12)}\qquad \{x:z_j(x;\theta_n)=1\}\subseteq\left\{x:\max_{j=g+1,\ldots,k} D_j(x;\theta_n)\geq D_1(x;\theta_n)\right\}$$

for $j=g+1,\ldots,k$, where the right-hand side converges toward the empty set, when $n$ tends to $\infty$, due to (A.3) and (A.4).

We can now use (A.11) in order to get

$$\limsup_{n\to\infty} L(\theta_n,P)\leq\lim_{n\to\infty} E_P\left[\sum_{j=1}^{g} z_j(\cdot;\theta_n)\left(\log\pi_j^n-\frac{p}{2}\log 2\pi-\frac{1}{2}\log|\Sigma_j^n|\right.\right.$$
$$\left.\left.-\frac{1}{2}(\cdot-\mu_j^n)'(\Sigma_j^n)^{-1}(\cdot-\mu_j^n)\right)\right]$$

$$= E_P\left[\sum_{j=1}^{g} z_j(\cdot; \tilde{\theta})\left(\log \pi_j - \frac{p}{2}\log 2\pi - \frac{1}{2}\log |\Sigma_j|\right.\right.$$

$$\left.\left. - \frac{1}{2}(\cdot - \mu_j)'\Sigma_j^{-1}(\cdot - \mu_j)\right)\right],$$

where $x \mapsto z_j(x; \tilde{\theta})$ are the assignment functions which would be derived when working with $g$ (instead of $k$) populations and $\tilde{\theta}$ being equal to a limit of the subsequence $\{\tilde{\theta}_n\}_{n=1}^{\infty} = \{(\pi_1^n, \ldots, \pi_g^n, \mu_1^n, \ldots, \mu_g^n, \Sigma_1^n, \ldots, \Sigma_g^n)\}_{n=1}^{\infty}$.

As $\sum_{j=1}^{g} \pi_j < 1$, the proof ends up by showing that we can change the weights $\pi_1, \ldots, \pi_k$ by

$$(A.13) \qquad \pi_j^* = \frac{\pi_j}{\sum_{j=1}^{g} \pi_j} \qquad \text{for } 1 \leq i \leq g \text{ and } \pi_{g+1}^* = \cdots = \pi_k^* = 0$$

(and properly modifying the assignment functions $z_j$'s). This change produces a strict decrease in the objective function, leading to a contradiction with the optimality stated in (2.7). Thus, we conclude $g = k$. $\square$

PROOF OF PROPOSITION 2.  Taking into account previous lemmas, we have that one of the two possibilities must hold.

(i) If $\pi_j^n \to \pi_j > 0$ for $1 \leq j \leq k$, then the choice of $\theta$ is obvious.

(ii) If $\pi_j^n \to \pi_j > 0$ with $\pi_j > 0$ for $j \leq g$ and $\pi_j = 0$ for $g < j \leq k$, we can define weights $\pi_j$ as $\pi_j = \lim_{n\to\infty} \pi_j^n$ for $j = 1, \ldots, g$ and $\pi_{g+1} = \cdots = \pi_k = 0$, and, take $\mu_j = \lim_{n\to\infty} \mu_j^n$ and $\Sigma_j = \lim_{n\to\infty} \Sigma_j^n$ for $j \leq g$. The other $\mu_j$'s and $\Sigma_j$'s may be arbitrarily chosen (of course, satisfying the eigenvalues-ratio restrictions). $\square$

**A.2. Consistency.** Given $\{x_n\}_{n=1}^{\infty}$ an i.i.d. random sample from an underlying (unknown) probability distribution $P$, let $\{\theta_n\}_{n=1}^{\infty} = \{(\pi_1^n, \ldots, \pi_k^n, \mu_1^n, \ldots, \mu_k^n, \Sigma_1^n, \ldots, \Sigma_k^n)\}_{n=1}^{\infty} \subset \Theta_c$ denote a sequence of sample estimators obtained by solving the problem (2.4) for the empirical measures $\{P_n\}_{n=1}^{\infty}$ with the eigenvalue-ratio restrictions (ER) for a fixed $c \geq 1$ [from Proposition 2 such a sequence does always exist, for large enough $n$, whenever $P$ is an absolutely continuous distribution verifying (PR)]. Notice that although similar notation to that applied in Section A.1 is here used, the index $n$ will now indicate the dependence on the sample size $n$.

The proof of the consistency combines arguments already used to prove the existence and techniques in the modern theory of empirical processes (see, e.g., Van der Vaart and Wellner [28]). We limit ourselves to state the key results as lemmas. The complete proofs can be obtained in [16].

First, we prove that there exists a compact set $K \subset \Theta_c$ such that $\theta_n \in K$ for $n$ large enough with probability 1. This follows from the next lemmas.

LEMMA A.4. *If $P$ is an absolutely continuous distribution [thus verifying (PR)] then the minimum (resp. maximum) eigenvalue, $m_n$ (resp. $M_n$) of the matrices $\Sigma_j^n$, $j = 1, \ldots, k$, cannot verify $m_n \to 0$ (resp. $M_n \to \infty$).*

LEMMA A.5. *If $P$ is an absolutely continuous distribution, then we can choose empirical centers $\mu_j^n$, $j = 1, \ldots, k$, such that their norms are uniformly bounded with probability 1.*

The following lemmas, related to uniform convergences, complete our technical needs for the final proof. The second involves $R(\theta; P)$ in (2.6).

LEMMA A.6. *Given a compact set $K$, the class of functions*

$$(A.14) \qquad \mathcal{H} := \left\{ I_{[u,\infty)}(D(\cdot; \theta)) \sum_{j=1}^{k} z_j^*(\cdot; \theta) \log D_j(\cdot; \theta) : \theta \in K, u \geq 0 \right\}$$

*is a Glivenko–Cantelli class, with $z_j^*(x; \theta) = I\{x : D(x; \theta) = D_j(x; \theta)\}$. (All the points in $\mathbb{R}^p$ are assigned to some class through the $z_j^*$'s.)*

LEMMA A.7. *Let $P$ be an absolutely continuous distribution with an strictly positive density function. Then for every compact subset $K$, we have that*

$$(A.15) \qquad \sup_{\theta \in K} |R(\theta; P_n) - R(\theta; P)| \to 0, \qquad P\text{-a.e.}$$

We can now prove the stated consistency result.

PROOF OF PROPOSITION 3. Let $K$ be a compact set such that $\theta_n \in K$ for $n \geq n_0$ with probability 1. The objective function in the empirical case can be rewritten as:

$$L(\theta, P_n) = \int_{\{x : D(x, \theta) \geq R(\theta; P_n)\}} \left[ \sum_{j=1}^{k} z_j^*(x; \theta) \log D_j(x; \theta) \right] dP_n(x),$$

with the $z_j^*$'s introduced in Lemma A.6. Let us introduce

$$\widetilde{L}(\theta, P_n) = \int_{\{x : D(x, \theta) \geq R(\theta; P)\}} \left[ \sum_{j=1}^{k} z_j^*(x; \theta) \log D_j(x; \theta) \right] dP_n(x).$$

We can see that

$$\sup_{\theta \in K} |L(\theta; P_n) - \widetilde{L}(\theta; P_n)| = o_P(1),$$

by using Lemma A.7 and the fact that the integrand can be bounded from above and below from some constants uniformly for $\theta$ in the compact set $K$.

Finally, we resort to the Glivenko–Cantelli property for the class of functions $\mathcal{H}$ in (A.14), and apply Theorem 3.2.3 in [28] to achieve the result. $\square$

**Acknowledgments.** We thank the associated editor and two anonymous referees for their valuable comments, which have greatly improved the paper.

DEPARTAMENTO DE ESTADÍSTICA E
INVESTIGACIÓN OPERATIVA
UNIVERSIDAD DE VALLADOLID
C/PRADO DE LA MAGDALENA, S/N
VALLADOLID 47005
SPAIN
E-MAIL: lagarcia@eio.uva.es
        alfonsog@eio.uva.es
        matran@eio.uva.es
        agustinm@eio.uva.es